\newcommand{\ie}{{\it i.e.\/}\ }
\newcommand{\ext}[1]{\stackrel{#1}{\wedge}}

\newcommand{\ctX}{\mbox{${\tilde{\cal X}}$}}

\def\limproj{\mathop{\oalign{lim\cr
\hidewidth$\longleftarrow$\hidewidth\cr}}}

\def\QQ{\rlap {\raise 0.4ex \hbox{$\scriptscriptstyle |$}}
  {\hskip -0.1em Q}}

\def\CC{{\rm \kern.24em \vrule width.02em height1.4ex
depth-.05ex \kern-.26em C}}
\def\C*{{\CC}^{*}}
\def\CQ{{{\CC}^{*}} \otimes {\QQ}}
\def\CM{{CM_{\infty}}}
\def\CMX{{CM_{\infty}(X)}}
\def\Z{{\bf Z}}

\def\wt{\widetilde}
\def\ot{\otimes }
\def\ra{\rightarrow }
\def\ify{\infty}

\def\Tg{{\cal T}_g}
\def\Tt{{\cal T}_{\tilde g}}
\def\Tbg{{\cal T}_{\bar g}}

\def\Tin{{\cal T}_{\infty}}
\def\TinX{{\cal T}_{\infty}(X)}
\def\TinY{{\cal T}_{\infty}(Y)}
\def\Cg{{\cal C}_g}

\def\Hin{{H}_{\infty}}

\def\Mg{{\cal M}_g}

\def\Dng{{DET}_{n,g}}
\def\Dnt{{DET}_{n,\tilde g}}
\def\Dnb{{DET}_{n,\bar g}}
\def\Dngt{{DET}_{n,g}^{\otimes 12}}
\def\Dntt{{DET}_{n,\tilde g}^{\otimes 12}}

\def\DnQ{{DET}(n,{\QQ})}
\def\D0Q{{DET}(0,{\QQ})}

\def\Tpi{{\cal T}(\pi)}

\def\tg{{\tilde g}}
\def\tX{{\tilde X}}
\def\tY{{\tilde Y}}

\def\bg{{\bar g}}
\def\bX{{\bar X}}

\def\G{{\cal G}}
\def\L{{\cal L}}
\def\M{{\cal M}}
\def\O{{\cal O}}

\def\S{{\cal S}}
\def\T{{\cal T}}
\def\X{{\cal X}}

\def\a{\alpha }
\def\b{\beta }
\def\ga{\gamma }
\def\Ga{\Gamma }

\def\r{\rho }
\def\p{\psi }

\def\l{\lambda }
\def\la{\lambda }
\def\La{\Lambda }
\def\a{\alpha }
\def\b{\beta }
\def\r{\rho }
\def\p{\psi }
\def\de{\delta }

\def\om{\omega }
\def\Om{\Omega }

\def\Tpi{{\cal T}(\pi)}
\long\def\comment#1\endcomment{}

\def\mapright#1{\smash{
    \mathop{\longrightarrow}\limits^{#1}}}

\def\mapdown#1{\Big\downarrow
 \rlap{$\vcenter{\hbox{$\scriptstyle#1$}}$}}
\documentstyle[12pt]{article}
\begin{document}
\baselineskip16pt
\begin{flushright}
{\underline {{\it To appear in:} {\bf Acta Mathematica}}}
\end{flushright}

\begin{center}
{\bf{DETERMINANT BUNDLES, QUILLEN METRICS, AND MUMFORD ISOMORPHISMS
OVER THE UNIVERSAL COMMENSURABILITY TEICHM\"ULLER SPACE}}\footnote
{{\it IHES/M/95/43 (Paris) preprint}, May 1995.}\\
\vspace{0.4cm}
{\bf Indranil Biswas, Subhashis Nag and Dennis Sullivan} \\
\end{center}
\begin{center}
{\bf Abstract}
\end{center}

A coherent, genus-independent, equivariant construction of determinant
line bundles, and connecting Mumford isomorphisms, is obtained over
the inductive limit of the Teichm\"uller spaces of Riemann surfaces of
varying genus. The direct limit of the Teichm\"uller spaces
corresponds to the inverse limit over the directed system of all
finite unbranched pointed coverings of a fixed surface.

\section{Introduction}

Let $\Tg$ denote the Teichm\"uller space comprising compact marked
Riemann surfaces of genus $g$. Let $DET_n \longrightarrow \Tg$ be
the determinant (of cohomology) line bundle on $\Tg$ arising from
the $n$-th tensor power of the relative cotangent bundle on the
universal family $\Cg$ over $\Tg$. The bundle $DET_0$ is called the Hodge line
bundle. The bundle $DET_n$ is equipped with a hermitian structure which is
obtained from the construction of Quillen of metrics on
determinant bundles using the Poincar\'e metric on the relative tangent
bundle of $\Cg$, [Q].

These natural line bundles over $\Tg$ carry liftings of the standard
action of the mapping class group, $MC_g$, on $\Tg$. We shall think
of them as $MC_g$-equivariant line bundles, and the isomorphisms
we talk about will be $MC_g$-equivariant isomorphisms.
By applying the Grothendieck-Riemann-Roch theorem,
Mumford [Mum] had shown that $DET_n$ is a certain fixed
(genus-independent) tensor power of the Hodge bundle. More precisely,
$$
DET_{n} ~=~ DET^{\otimes (6n^2-6n+1)}_{0}
\eqno{(1.1)}
$$
this isomorphism of equivariant line bundles being ambiguous only
up to multiplication by a non-zero scalar. (Any choice of such an
isomorphism will be called a Mumford isomorphism in what follows.)

There is a remarkable connection, discovered by Belavin and Knizhnik [BK],
between the Mumford isomorphism above for the case $n=2$, [i.e., that
$DET_{2}$ is the $13$-th tensor power of Hodge], and the existence of the
Polyakov string measure on the moduli space $\Mg$. (See the discussion after
Theorem 5.5 for more details.) This suggests the question of finding a
genus-independent formulation of the Mumford isomorphisms over some
``universal'' parameter space of Riemann surfaces (of varying topology).

In this paper we combine a Grothendieck-Riemann-Roch lemma (Theorem 2.9)
with a new concept of $\CQ$ bundles (Section 5), to construct a
universal version of the determinant bundles and Mumford's isomorphism.
Our objects exist over a universal base space $\Tin = \TinX$, which is the
infinite directed union of the complex manifolds that are the Teichm\"uller
spaces
of higher genus surfaces that are unbranched coverings of any (pointed)
reference surface $X$. The bundles and the relating isomorphisms
are equivariant with respect to the natural action of the
universal commensurability group $\CM$ -- which is defined (up to isomorphism)
as the group of isotopy classes of unbranched self correspondences of the
surface $X$ arising from pairs of non-isotopic pointed covering maps
$X' {\ra \atop \ra} X$, (see below and in Section 5).

In more detail, our universal objects are obtained by taking the direct
limits using the following category $\S$: for each integer $g \geq 2$,
there is one object in $\S$, an oriented closed pointed surface $X_g$
of genus $g$, and one morphism $X_{\tilde g} \ra X_g$ for each based
isotopy class of finite unbranched pointed covering map. For each morphism
of $\S$ (say of degree $d$) we have an induced holomorphic injection of
Teichm\"uller spaces arising from pullback of complex structure:
$$
\Tpi:~\Tg \ra \Tt
\eqno(1.2)
$$
The GRR lemma provides a natural isomorphism of the line bundle
$DET_{n,g}^{\otimes 12d}$ on $\Tg$ with the pullback line bundle
${\Tpi}^*DET_{n,\tg}^{\otimes 12}$. We may view this isomorphism,
equivalently, as a degree $d$ homomorphism covering the injection $\Tpi$
between the principal $\C*$ bundles associated to $DET_{n,g}$ and
$DET_{n,\tg}$, respectively.  Each commutative triangle in $\S$ yields
a commutative triangular prism whose top face is the following triangle
of total spaces of principal $\C*$ bundles:
$$
\matrix{
\Dng^{\otimes {12}}
&
&\mapright{}
&
&\Dnt^{\otimes {12}}
\cr
&
\searrow
&
&\swarrow
&
\cr
&&
\Dnb^{\otimes {12}}
&&
\cr}
$$
and whose bottom face is the commuting triangle of base spaces for these
bundles:
$$
\matrix{
\Tg~~~~~
&
&\mapright{}
&
&~~~~~\Tt
\cr
&
\searrow
&
&\swarrow
&
\cr
&&
\Tbg
&&
\cr}
\eqno(1.3)
$$
Moreover, the canonical mappings above relating these $DET_n$ bundles
over the various Teichm\"uller spaces preserve the Quillen
hermitian structure of these bundles in the sense that unit circles are
carried to unit circles.

We explain in brief the commensurability Teichm\"uller space $\Tin$ and the
large mapping class group $\CM$ acting thereon. For
each object $X$ in $\S$, consider the directed set
$\{\a\}$
of all morphisms in $\S$ with range $X$. Then we form
$$
\TinX := dir. lim. \T_{g(\a)}
\eqno(1.4)
$$
where the limit is taken over $\{\a\}$, and $g(\a)$ is the genus for
the domain of morphism $\a$. Each morphism
$X_{g'} \ra X_g$ induces a holomorphic {\it bijection} of the
corresponding direct limits, and we denote any of these isomorphic
``ind-spaces'' (inductive limit of finite-dimensional complex spaces
-- see [Sha]) by $\Tin$ -- the universal commensurability
Teichm\"uller space. (Compare Section 2 and Example 4 on p. 547 of [S].)
Notice that a pair of morphisms $X' {\ra \atop \ra} X$
determines an automorphism of $\Tin$; we call the group of automorphisms
of $\Tin$ obtained this way the {\it commensurability modular group} $\CM$.

We now take the direct limit of the $\C*$ principal bundles associated
to $\Dngt$ over $\Tg$ to obtain a new object -- a $\CQ$ bundle
over $\Tin$ -- denoted $\DnQ$. As sets the total space with action of
the group $\CQ$ is defined by the direct limit construction. Continuity
and complex analyticity for maps into these sets are defined by the
corresponding properties for factorings through the strata of
the direct system. (Section 5.)

There are the natural isomorphisms of Mumford, as stated in (1.1),
at the finite-dimensional stratum levels. By our construction
these isomorphisms are {\it rigidified} to be natural over the
category $\S$. Therefore we have natural Mumford isomorphisms
between the following $\CQ$ bundles over the universal commensurability
Teichm\"uller space $\Tin$:
$$
\DnQ~~ {\rm and}~~ {\D0Q}^{\otimes (6n^2-6n+1)}
\eqno(1.5)
$$
We also show that the natural Quillen metrics of the $DET$ bundles
fit together to define a natural analogue of Hermitian structure on
these $\CQ$ bundles; in fact, for all our canonical mappings in
the direct system the unit circles are preserved. Note Theorem 5.5.

In fact, the existence of the canonical relating morphism between
determinant bundles (fixed $n$) in the fixed covering
$\pi:X_{\tilde {g}} \ra X_g$
situation was first conjectured and deduced by us utilizing the
differential geometry of these Quillen metrics. Recall that the
Teichm\"uller spaces $\Tg$ and $\Tt$ carry natural symplectic
forms (defined using the Poincare metrics on the Riemann surfaces)
-- the Weil-Petersson K\"ahler forms -- which are in fact
the curvature forms of the natural Quillen metrics of these $DET$
bundles ([Wol], [ZT], [BGS]). If the covering $\pi$ is unbranched of
degree $d$, a direct calculation shows that this natural WP
form  on $\Tt$ (appropriately renormalized by the degree)
pulls back to the WP form of $\Tg$ by $\Tpi$ (the embedding of
Teichm\"uller spaces induced by $\pi$).
One expects therefore that if one raises the $DET_n$ bundle on $\Tg$
by the tensor power $d$, then it extends over the larger Teichm\"uller
space $\Tt$ as the $DET_n$ bundle thereon. This intuition is, of course,
what is fundamentally behind our direct limit constructions. Since it
turns out to be technically somewhat difficult to actually prove that the
relevant bundles are isomorphic using this differential geometric
method, we have separated that aspect of our work into a different
article [BN].

Can objects on $\Tin$ that are equivariant by the commensurability modular
group $\CM$ be viewed as objects on the quotient ${\Tin}/{\CM}$ ?
This quotient is problematical and interesting, so we work with the
equivariant statement.

We end the Introduction by mentioning some problems.
The universal commensurability Teichm\"uller space,
$\T_{\ify}$, is made up from embeddings $\Tpi$ that are isometric
with respect to the natural
Teichm\"uller metrics, so it carries a natural Teichm\"uller metric.
Our theorems give us genus independent
determinant line bundles $\DnQ$, Quillen metrics and Mumford
isomorphisms over $\Tin$, all compatible with
each other and the commensurability group $\CM$.
{\it Are the above structures uniformly continuous for this metric?}
Then they would pass to the completion $\wt{T}_{\ify}$ of $\Tin$ for the
Teichm\"uller metric. One knows that $\wt{T}_{\ify}$ is a separable complex
Banach manifold which is the Teichm\"uller space of complex structures
on the universal solenoidal surface $H_{\infty}=\limproj {\tilde X}$,
where ${\tilde X}$ ranges, as above, over all finite covering
surfaces of $X$. (See [S], [NS], for the Teichm\"uller theory of
$H_{\infty}$.) We
would conjecture that this continuity is true and that the $\CQ$
bundles $\DnQ$, Quillen metrics, and Mumford isomorphisms can be
defined over $\T(H_{\infty})=\wt{T}_{\ify}$ directly by looking
at compact solenoidal Riemann surfaces themselves. Now we may consider
square integrable holomorphic forms on the leaves of $\Hin$ (which are
uniformly distributed copies of the hyperbolic plane in $\Hin$)
regarded as modules over the $\CC*$ algebra of $\Hin$ with chosen
transversal. The measure of this transversal would become a real
parameter extending the genus above. One expects that A. Connes'
version of Grothendieck-Riemann-Roch would replace Deligne's functorial
version which we are using here.

Finally one would hope that
the Polyakov measure (Section 5) on Teichm\"uller space, when viewed as a
metric on the canonical bundle, would also make sense at infinity in
the direct limit because this measure can be constructed by applying
the $13$-th power Mumford isomorphism ((1.1) for $n=2$)
to the $L^2$ inner product on the Hodge line bundle.
That issue remains open.

\noindent {\it Acknowledgements:} We would like express our gratitude
to several mathematicians
for their interest and discussions. In particular, we thank
M.S.Narasimhan and E.Looijenga for helpful discussions in early
stages of this work. Laurent Moret-Bailly deserves a very special
place in our acknowledgements. In fact, he brought to our attention
the Deligne pairing version of the Grothendieck-Riemann-Roch theorem
that we use crucially here, and showed us Lemma 2.9,
after seeing an earlier (less strong version) of the main
theorem based on topology and the curvature calculations mentioned above.


\section{A lemma on determinant bundles}

Let $X$ be a compact Riemann surface, equivalently, an
irreducible smooth projective curve over $\CC$. Let $L$ be a
holomorphic line bundle on $X$.  The determinant of $L$ is the
defined to be the $1$-dimensional complex vector space
$(\ext{top}H^0(X,L))\bigotimes (\ext{top}H^1(X,L)^*)$, and will
be denoted by $det(L)$. Take a Riemannian metric $g$ on $X$
compatible with the conformal structure of $X$. Fix a hermitian
metric $h$ on $L$. Using $g$ and $h$, a hermitian structure can
be constructed on ${\Om}^i(X,L)$, the space of $i$-forms on $X$
with values in $L$.  Moreover the vector space $H^1(X,L)$ is
isomorphic, in a natural way, with the space of harmonic
$1$-forms with values in $L$. Consequently the vector spaces
$H^0(X,L)$ and $H^1(X,L)$ are equipped with hermitian structures
which in turn induce a hermitian structure on $det(L)$ -- this
metric on $det(L)$ is usually called the $L^2$ metric. Let
$\Delta := {\overline\partial}^*{\overline\partial}$ be the
laplacian acting on the space of smooth sections of $L$. Let
$\{{\la}_i\}_{i\geq 1}$ be the set of non-zero eigenvalues of
$\Delta$; let $\zeta$ denote the analytic continuation of the
function $s \longmapsto \sum_{i}1/{\la}^s_i$. The {\it Quillen
metric} on $det(L)$ is defined to be the hermitian structure on
$det(L)$ obtained by multiplying the $L^2$ metric with ${\rm
exp}(-{\zeta}'(0))$, [Q].

To better suit our purposes, we will modify the above (usual)
definition of the Quillen metric by a certain factor.
Consider the real number $a(X)$ appearing in Th\'eor\`eme 11.4 of [D].
This number $a(X)$ depends only on the genus of $X$. The statement in
Remark 11.5 of [D] -- to the effect that there is a constant $c$ such
that $a(X) = c.\chi (X)$, where $\chi (X)$ is the Euler characteristic
of $X$ --  has been established in [GS]. (The constant $c$ is related
to the derivative at $(-1)$ of the zeta function for the trivial
hermitian line bundle on $\CC P^1$ (4.1.7 of [GS]).) Let $H_Q(L; g,h)$
denote the Quillen metric on $det(L)$ defined above. Henceforth, by
{\it Quillen metric} on $det(L)$ we will mean the hermitian metric
$$
{\rm exp}(a(X)/12) H_Q(L; g, h)
\eqno{(2.1)}
$$

Next we will describe briefly some key properties of the
determinant line and the Quillen metric.

Let $\pi:\X \longrightarrow S$ be a family of compact Riemann
surfaces parametrized by a base $S$. We can work with either
holomorphic (Kodaira-Spencer) families over a complex-analytic
variety $S$, or with algebraic families over complex algebraic
varieties (or, more generally, over a scheme) $S$. In the
algebraic category one means that $\pi$ is a proper smooth
morphism of relative dimension one with geometrically connected
fibers. In the analytic category, $\pi$ is a holomorphic
submersion again with compact and connected fibers.  Take a
hermitian line bundle $L_S \longrightarrow \X$ with hermitian
metric $h_S$. Fix a hermitian metric $g_S$ on the relative
tangent bundle $T_{\X/S}$.

For any point $s \in S$, the above construction gives a
hermitian line $det(L_s)$ (the hermitian structure is given by
the Quillen metric). The basic fact is that these lines fit
together to give a line bundle on $S$ [KM], which is called the
{\it determinant bundle} of $L_S$, and is denoted by $det(L_S)$.
Moreover the function on the total space of $det(L_S)$ given by
the norm with respect to the Quillen metric on each fiber is a
$C^{\infty}$ function, and hence it induces a hermitian metric
on $det(L_S)$ [Q]. This bundle will be denoted by $det(L_S)$.
We shall make clear in Remark 2.13 below that this ``determinant
of cohomology'' line bundle is also an algebraic or analytic
bundle -- according to the category within which we work.

The determinant bundle $det(L_S)$ is functorial with respect to
base change. We describe what this means. For a morphism $\ga :
S' \longrightarrow S$ consider the bundle, $p^*_2L_S
\longrightarrow S'\times_S\X$, on the fiber product, where $p_2
: S'\times_S\X \longrightarrow \X$ is the projection onto the
second factor. The hermitian structure $h_S$ pulls back to a
hermitian structure on $L_{S'} := p^*_2L_S$; and, similarly, the
metric $g_S$ induces a hermitian structure on the relative
tangent bundle of $S'\times_S\X$. ``Functorial with respect to
base change'' now means that in the above situation there is a
{\it canonical isometric isomorphism} $$ {\rho}_{S',S} ~:~
det(L_{S'})~\longrightarrow ~ {\ga}^*det(L_S)$$ such that if
$$S''~\mapright{\ga'}~S'~\mapright{{\ga}}~S$$ are two morphisms
then the following diagram is commutative $$
\matrix{det(L_{S''})&\mapright{{\rho}_{S'',S'}}&
{\ga'}^*(det(L_{S'}))&\cr \mapdown{{\rho}_{S'',S}}&&
\mapdown{{\ga'}^*{\rho}_{S',S}}&\cr (\ga\circ{\ga}')^*det(L_S)
&\mapright{id}& (\ga\circ{\ga}')^*det(L_S)&\cr}
\eqno(2.2)
$$

The determinant of cohomology construction $det(L_S)$ produces a
bundle over the parameter space $S$ induced by the bundle over
the total space $\X$; now, the Grothendieck-Riemann-Roch (GRR)
theorem gives a canonical isomorphism of $det(L_S)$ with a
combination of certain bundles obtained (on $S$) from the direct
images of the bundle $L_S$ and the relative tangent bundle
$T_{\X/S}$. In order to relate canonically the determinant
bundle obtained from a given family $\X \rightarrow S$ (fibers
of genus $g$, say) with the determinant arising from a covering
family $\ctX$ (having fibers of some higher genus $\tg$), we
shall utilize the GRR theorem in a formulation due to Deligne,
[D, Theorem 9.9(iii)].

In fact, Deligne introduces a ``bilinear pairing'' that
associates a line bundle, denoted by $<L_S,M_S>$, over $S$ from
any pair of line bundles $L_S$ and $M_S$ over the total space of
the fibration $\X \rightarrow S$. If $L_S$ and $M_S$ carry
hermitian metrics then a canonically determined hermitian
structure gets induced on the Deligne pairing bundle $<L_S,M_S>$
as well. Denoting by $\L = L_S$ the given line bundle over $\X$,
the GRR theorem in Deligne's formulation reads:


$$ det(\L)^{\otimes 12}~=~ <T^*_{\X /S},
T^*_{\X /S}> \bigotimes <{\L} , \L {\otimes} T_{\X /S}>^{\otimes
6}
\eqno(2.3)
$$
Here $T^*_{\X /S}$ denotes the
relative cotangent bundle over $\X$, and the equality asserts
that there is a {\it canonical isomorphism, functorial with
respect to base change,} between the bundles on the two sides.
Furthermore, Th\'eor\`eme 11.4 of [D] says that the canonical
identification in (2.3) is actually an {\it isometry} with the
Quillen metric on the left side and the Deligne pairing metrics
on the right. (The constant ${\rm exp}(a(X))$ in the statement
of Th\'eor\`eme 11.4 of [D] has been absorbed in the definition
(2.1).) We proceed to explain the Deligne pairing and the metric
thereon in brief; details are to be found in sections 1.4 and
1.5 of [D].

Let $L$ and $M$ be two line bundles on a compact Riemann surface
$X$. For a pair of meromorphic sections $l$ and $m$ of $L$ and
$M$ respectively, with the divisor of $l$ being disjoint from
the divisor of $m$, let $\CC <l,m>$ be the one dimensional
vector space with the symbol $<l,m>$ as the generator. For two
meromorphic functions $f$ and $g$ on $X$ such that $div(f)$ is
disjoint from $div(m)$ and $div(g)$ is disjoint from $div(l)$,
the following identifications of complex lines are to be made
$$
\matrix{<fl,m> &= & f(div(m))<l,m>&\cr
<l,gm>&=&g(div(l))<l,m>&}
\eqno{(2.4)}
$$ The Weil reciprocity law says that for any two meromorphic
functions $f_1$ and $f_2$ on $X$ with disjoint divisors,
$f_1(div(f_2)) = f_2(div(f_1))$ [GH, page 242]. So we have
$$<fl,gm>~=~f(div(gm)).g(div(l))<l,m>~=~g(div(fl)).f(div(m))<l,m>.$$
{}From the above equality it follows that the identifications in
(2.4) produce a complex one dimensional vector space, denoted by
$<L,M>$, from the pair of line bundles $L$ and $M$. If $L$ and
$M$ are both equipped with hermitian metrics then the hermitian
metric on $\CC <l,m>$ defined by
$$
log||<l,m>|| := {1\over{2{\pi}i}}\int_X\partial
{\overline\partial}(log||l||.log||m||)~+\, log||l||(div(m))
\,+\,log||m||(div(l))
\eqno{(2.5)}
$$
is compatible with the relations in (2.4) --  hence it gives
a hermitian structure on $<L,M>$, see [D, 1.5.1].

Consider now a family of Riemann surfaces $\X \longrightarrow
S$; let $L_S$ and $M_S$ be two line bundles on $\X$, equipped
with hermitian structures. Over an open subset $U\subset S$, let
$l_U$, $m_U$ be two meromorphic sections of $L_S$ and $M_S$
respectively, with finite supports over $U$ such that the
support of $l_U$ is disjoint from the support of $m_U$. (Support
of a section is the divisor of the section.) For another open
set $V$ and two such sections $l_V$ and $m_V$, the relations in
(2.4) give a function $$C_{U,V}~\in ~{\O}^*_{U\cap V}\,.$$ Using
the Weil reciprocity law it can be shown that $\{C_{U,V}\}$
forms a 1-cocycle on $S$. In other words, we get a line bundle
on $S$, which we will denote by $<L_S,M_S>$. The hermitian
structure on $<L,M>$, described earlier, makes $<L_S,M_S>$ into
a hermitian bundle.

Given a meromorphic section $m$ of $M_S$, let $m^{\otimes n}$ be
the meromorphic section of $M^n_S$ obtained by taking the $n$-th
tensor power of $m$. Note that $div(m^{\otimes n}) = n.div(m)$.
The map $<l,m^{\otimes n}> \longrightarrow <l,m>^{\otimes n}$
can be checked to be compatible with the relations (2.4) and
hence it induces an isomorphism
$$
<L_S,M^n_S> \longrightarrow <L_S,M_S>^n
\eqno{(2.6)}
$$ From the definition (2.5) we see that (2.6) is an isometry
for the metric on $M^n_S$ induced by the metric on $M_S$.

We shall now see how the critical formula (2.3) follows from the
general GRR theorem of [D]. Indeed, let $\L$ denote any rank $n$
vector bundle on the total space of the family $\X$; we
reproduce below the statement of Theorem 9.9(iii) of [D]:
$$
det(\L )^{\otimes 12}~=~ <T^*_{\X /S}, T^*_{\X /S}> \bigotimes
<{\Lambda^{n}}(\L ),{\Lambda^{n}}(\L ) {\otimes} T_{\X
/S}>^{\otimes 6} \bigotimes I_{\X /S}C^2(\L )^{-12}
\eqno{(GRR-D)}
$$
Now, from the definition of $I_{\X /S}C^2$ in [D, 9.7.2] it
follows that for a line bundle $\L$, the bundle
$I_{\X /S}C^2(\L)$ is the
trivial bundle on $S$, and the metric on it is the constant
metric [D, Theorem 10.2(i)]. From Th\'eor\`eme 11.4 of [D] we
conclude that that the canonical identification in the statement
above is actually an isometric identification.  (The factor
${\rm exp}(a(X))$ in Th\'eor\`eme 11.4 of [D]
is taken care of by the definition (2.1).) Thus we have obtained
the isometric isomorphism stated in (2.3).

With this background behind us we can formulate our main lemma.
Let $\X$ and $\ctX$ be two families of Riemann surfaces over $S$
(say with fibers of genus $g$ and $\tg$, respectively), and $p:
\ctX \longrightarrow \X$ be an \'etale (\ie unramified) covering
of degree $d$, commuting with the projections onto $S$. In other
words, the map $p$ fits into the following commutative diagram
$$
\matrix{\ctX& &\mapright{p}& & \X
\cr &\searrow& & \swarrow &
\cr & &S& & \cr}
\eqno(2.7)
$$ The situation implies that each fiber of the family $\ctX$ is
a degree $d=(\tg - 1)/(g - 1)$ holomorphic covering over the
corresponding fiber of the family $\X$.  Fix also a hermitian
metric $g$ on $T_{\X /S}$.  Since $p$ is \'etale, $p^*T_{\X /S}
= T_{\ctX /S}$, and hence $g$ induces a hermitian metric $p^*g$
on $T_{{\ctX}/S}$. Let ${\X}' \longrightarrow S$ be a third
family of Riemann surface which is again an \'etale cover of
$\ctX$ and fits into the following commutative diagram $$
\matrix
{{\X}' &\mapright{q} &\ctX &\mapright{p}& \X
\cr &\searrow& \mapdown{} &\swarrow &
\cr & &S && \cr}
\eqno(2.8)
$$ We want to prove the following:


\medskip
\noindent{\bf Lemma 2.9.}\,\ {\it (i)~Let $\L$ be a hermitian
line bundle on $\X$ and let $p^*\L \longrightarrow \ctX$ be the
pullback of $\L$ equipped with the pullback metric. Then there
is a canonical isometric isomorphism $$det(p^*(\L ))^{\otimes
12}~\cong ~ det(\L )^{\otimes 12.{\rm deg}(p)}$$ of bundles on
$S$. This isomorphism is functorial with respect to base change.

\noindent (ii)~ Denoting the isometric isomorphism obtained in
(i) by ${\Ga}(p)$, and similarly defining ${\Ga}(q)$ and
${\Ga}(p\circ q)$, the following diagram commutes:
$$\matrix{det((p\circ q)^*(\L ))^{\otimes
12}&\mapright{{\Ga}(q)}& det(p^*(\L ))^{\otimes 12.{\rm
deg}(q)}\cr \mapdown{{\Ga}(p\circ q)}&&
\mapdown{{\Ga}(p)^{\otimes {\rm deg}(q) }}\cr det(\L )^{\otimes
12.{\rm deg}(p\circ q)}&\mapright{id}& det(\L )^{\otimes 12.{\rm
deg}(p\circ q)}\cr} $$ where ${\Ga}(p)^{\otimes {\rm deg}(q) }$
is the isomorphism on appropriate bundles, obtained by taking
the ${\rm deg }(q)$-th tensor product of the isomorphism
${\Ga}(p)$.}
\medskip

(The terminology ``functorial with respect to base change'' was
explained earlier. We will use ``canonical'' to mean functorial
with respect to base change.)

\medskip
\noindent {\bf Proof of Lemma 2.9.} The idea of the proof is to
relate -- utilizing GRR in form (2.3) -- the determinant
bundles, which are difficult to understand, with the more
tractable ``Deligne pairings''.

Let $\M$ be any line bundle on $\ctX$ equipped with a hermitian
structure. First we want to show that there is a canonical
isometric isomorphism

$$<p^*\L ,\M >~\longrightarrow ~<\L ,N(\M )>\, , \eqno{(2.10)}$$
where $N(\M ) \longrightarrow \X$ is the norm of $\M$. We recall
the definition of $N(\M)$. The direct image $R^0p_*(\M )$ is
locally free on $\X$, and the bundle $R^0p_*(\M )$ admits a
natural reduction of structure group to the {\it monomial group}
$G \subset GL({\rm deg}(p), \CC)$. (The group $G$ is the
semi-direct product of permutation group, $P_{{\rm deg}(p)}$,
with the invertible diagonal matrices defined using the
permutation action of $P_{{\rm deg}(p)}$.) Mapping $g \in G$
to the permanent of $g$ (on $G$ it is simply the product
of all non-zero entries) we get a homomorphism to ${\CC}^*$,
which is denoted by $\mu$. Using this homomorphism $\mu$
we have a holomorphic line bundle on $\X$, associated to
$R^0p_*(\M )$, which is defined to be $N(\M)$. Clearly the
fiber of $N(\M)$ over a point $x\in \X$ is the tensor product
$$ N(\M )_x ~ = ~
\bigotimes_{y\in p^{-1}(x)}{\M}_y\,,
\eqno{(2.11)}
$$
The hermitian metric on $\M$ gives a reduction of
the structure group of $R^0p_*(\M )$ to the maximal compact
subgroup $G_U  \subset  G$. Since $\mu (G_U) = U(1)$, we have
a hermitian metric on $N(\M)$. Note that the hermitian metric on
$N(\M )$ is such that the above equality (2.11) is actually an
isometry.

For a meromorphic section $m$ of $M$, the above identification
of fibers gives a meromorphic section of $N(M)$ which is denoted
by $N(m)$. Given sections $l$ and $m$ of $\L$ and $\M$
respectively, with finite support over $U\subset S$ (the support
of $p^*l$ and $m$ being assumed disjoint) we map $<p^*l,m>$ to
$<l,N(m)>$. It can be checked that this map is compatible with
the relations in (2.4). Hence we get a homomorphism from the
bundle $<p^*\L ,\M>$ to $<\L ,N(\M )>$; this is our candidate
for (2.10). To check that it is an isometry, we evaluate the
(logarithms of) norms of the sections $<p^*l,m>$ and $<l,N(m)>$
given by definition (2.5). It is easy enough to see from (2.5) that
the norms of these two sections coincide.

Therefore for a hermitian line bundle ${\L}'$ on $\X$, the
isomorphism (2.10) implies that $$<p^*\L ,p^*{\L}'> ~=~ <\L
,N(p^*{\L}')>$$ But $N(p^*{\L}') = {\L'}^d$, where $d :={\rm
deg}(p)$, and moreover the hermitian metric on $N(p^*{\L}')$
coincides with that of ${\L'}^d$. Hence from the isometric
isomorphism obtained in (2.6) we get the following
identification of hermitian line bundles (the isomorphism so
created being again functorial with respect to change of base
space): $$ <p^*\L,p^*{\L}'> ~=~<\L ,{\L}'>^d
\eqno{(2.12)}
$$

To prove part $(i)$ of the Lemma we apply the GRR isomorphism
(2.3) to both $\L$ and $p^*\L$, and compare the Deligne pairing
bundles appearing on the right hand sides using the result
(2.12).  To simplify notation set $\om = T^*_{\X /S}$. By
applying (2.3) to $p^*\L$, and noting that since the map $p$ is
\'etale, the relative tangent bundle $T_{\ctX /S} = p^*T_{\X
/S}$, we deduce that $det(p^*\L )^{\otimes 12}$ is canonically
isometrically isomorphic to $<p^*\L ,p^*(\L\bigotimes{\om}^{-1}
)>^{\otimes 6}
\bigotimes <p^*\om ,p^*{\om}>$.
Taking ${\L}'$ to be ${\L}\bigotimes\om$ in (2.12) we have
$<\L,\L\bigotimes{\om}>^d=<p^*\L,(p^*\L\bigotimes \om)>$.
Substituting $\om$ in place on $\L$ and ${\L}'$ in (2.12) we
have $<\om ,\om>^d=<p^*\om ,p^*\om >$. Therefore the bundle
$<p^*\L ,p^*(\L\bigotimes{\om}^{-1})>^{\otimes 6}\bigotimes
<p^*\om ,p^*{\om}>$ is isometrically isomorphic to
$<\L,\L\bigotimes{\om}^{-1}>^{6d}\bigotimes <\om ,\om>^d$. But
now applying (2.3) to $\L$ itself we see that this last bundle
is isometrically isomorphic to $det(\L )^{\otimes 12d}$. That
completes the proof. Notice that since all isomorphisms used in
the above proof were canonical (functorial with base change),
the final isomorphism asserted in part $(i)$ is also canonical
in the same sense.

In order to prove part $(ii)$ of the Lemma, we first note that
the isometric isomorphisms in (2.10) and (2.12) actually fit
into the following commutative diagram $$
\matrix{<(p\circ q)^*\L ,\M >&\mapright{}&<p^*\L ,N(\M )_q>
\cr \mapdown{}&& \mapdown{}\cr <\L ,N({\M})>& = & <\L
,N({\M})>\cr} $$
where $\L$ is a hermitian line bundle on $\X$ and $\M$ is a
hermitian line bundle on ${\X}'$, $N(\M ) \longrightarrow \X$ is
the norm of $\M$ for the covering $p\circ q$, and $N(\M )_q
\longrightarrow \ctX$ is the norm of $\M$ for the covering $q$.
Indeed, the commutativity of the above diagram is
straightforward to deduce from the fact that the following two
bundles on $\X$: namely, $N(\M )$ and the norm of $N(\M )_q$,
are isometrically isomorphic. The isomorphism can be defined,
for example, using (2.11). Now using (2.3), and repeatedly using
the above commutative diagram, we obtain part $(ii)$.
$\hfill{\Box}$
\medskip

We will have occasion to use this general lemma in concrete
situations.

\medskip
\noindent{\bf Remark 2.13.} In [KM] and in [D] the basic context
is the algebraic families category, and the determinant of
cohomology bundle as well as the Deligne pairing bundles are
constructed in this category.  However, since the constructions
of the determinant bundles and of the Deligne pairing are {\it
canonical} and {\it local}, they work equally well for
holomorphic families of Riemann surfaces also. The point is that
if $\X \rightarrow S$ is a holomorphic family of Riemann
surfaces parametrized by a complex manifold $S$, and $\L
\rightarrow \X$ is a holomorphic line bundle, then $det(\L)
\rightarrow S$ is a holomorphic line bundle which is functorial
with respect to holomorphic base changes. And if $\L$ and $\M$
are two holomorphic line bundles on $\X$ then $<\L,\M>$ is a
holomorphic line bundle on $S$ -- again functorial with respect
to holomorphic base changes. In fact, an analytic construction of
the determinant bundle and the Quillen metric is to be found in [BGS].

Since the constructions of the Quillen metric and the metric on
the Deligne pairing, (using (2.5)), also hold true for
holomorphic families, consequently, Lemma 2.9 holds in the
holomorphic category as well as in the algebraic category.

\noindent{\bf Remark 2.14.} The statement that
$det(p^*(\L))^{\otimes 12}~\cong ~ det(\L )^{\otimes 12{\rm
deg}(p)}$ as line bundles actually holds for curves over any
field. The statement about isometry makes sense only when we
have Riemann surfaces.

\section{Determinant bundles over Teichm\"uller spaces}

Our aim in this section is to apply the Lemma 2.9 to the universal
family of marked Riemann surfaces of genus $g$ over the Teichm\"uller
space $\Tg$. The situation of Lemma 2.9 is precipitated by choosing
any finite covering space over a topological surface of genus $g$.

Let $\pi:\tX \longrightarrow X$ be an unramified covering map between
two compact connected oriented two manifolds $\tX$ and $X$ of genera
$\tg$ and $g$, respectively. Assume that $g \geq 2$. The degree of the
covering $\pi$, which will play an important role, is the ratio
of the respective Euler characteristics;
namely, $deg(\pi)=(\tg-1)/(g-1)$.

We recall the basic deformation spaces of complex (conformal)
structures on smooth closed oriented surfaces -- the Teichm\"uller spaces.
Let $Conf(X)$ (respectively, $Conf(\tX)$) denote the space of
all smooth conformal structures on $X$ (respectively, $\tX$). Define
${\rm Diff}^{+}(X)$ (respectively, ${\rm Diff}^+(\tX)$)
to be the group of all orientation preserving diffeomorphisms of
$X$ (respectively, $\tX$), and denote by
${\rm Diff}^{+}_0(X)$ (respectively, ${\rm Diff}^+_0(\tX)$)
the subgroup of those that are homotopic to the identity.

The group ${\rm Diff}^+(X)$ acts naturally on $Conf(X)$ by pullback
of conformal structure. We define
$$
{\cal T}(X)={\cal T}_g~:=~Conf(X)/{\rm Diff}^{+}_0(X)
\eqno{(3.1)}
$$
to be the Teichm\"uller space of genus $g$ (marked) Riemann surfaces.
Similarly obtain ${\cal T}_{\tg} := Conf(\tX )/{\rm Diff}^+_0(\tX )$
-- the Teichm\"uller space for genus $\tg$. The Teichm\"uller space
$\Tg$ carries naturally the structure of a $(3g-3)$ dimensional
complex manifold which is embeddable as a contractible domain
of holomorphy in the affine space $\CC^{3g-3}$. The mapping class
group of the genus $g$ surface, namely the discrete group
$MC_g := {\rm Diff}^+(X)/{\rm Diff}^+_0(X)$, acts properly
discontinuously on $\Tg$ by holomorphic automorphisms, the quotient
being the moduli space $\Mg$. For these basic facts see, for example,
[N].

The Teichm\"uller spaces are fine moduli spaces. In fact,
the total space $X\times {\cal T}_g$ admits a natural complex
structure such that the projection to the second factor
$$
{\psi}_g:{\cal C}_g:=X\times{\cal T}_g \longrightarrow {\cal T}_g
\eqno{(3.2)}
$$
gives the universal Riemann surface over ${\cal T}_g$. This means that for
any $\eta \in {\cal T}_g$, the submanifold $X\times\eta$ is a complex
submanifold of ${\cal C}_g$, and the complex structure on $X$ induced by
this embedding is represented by $\eta$. As is well-known, (Chapter 5
in [N]), the family $\Cg \rightarrow \Tg$ is the {\it universal} object
in the category of holomorphic families of genus $g$ marked
Riemann surfaces.

Given a complex structure on $X$, using $\pi$ we may pull back this to a
complex structure on $\tX$. This gives an injective map
$Conf(X)\longrightarrow Conf({\tX})$. Given an element of $f\in {\rm
Diff}^+_0(X)$, from the homotopy lifting property, there is a unique
diffeomorphism ${\tilde f}\in{\rm Diff}^{+}_0({\tX})$ such that $\tilde f$
is a lift of $f$. Mapping $f$ to $\tilde f$ defines an injective
homomorphism of ${\rm Diff}^{+}_0(X)$ into ${\rm Diff}^{+}_0({\tX})$.
We therefore obtain an injection
$$
{\cal T}(\pi ):{\cal T}_g~\longrightarrow ~{\cal T}_{\tg}
\eqno(3.3)
$$
It is known that this map ${\cal T}(\pi )$ is a {\it proper holomorphic
embedding} between these finite dimensional complex manifolds; $\Tpi$
respects the quasiconformal-distortion (=Teichm\"uller) metrics.
{}From the definitions it is evident that this embedding between the
Teichm\"uller spaces depends only on the (unbased) isotopy class of
the covering $\pi$.

\noindent{\bf Remark 3.4.} In fact, we see that
$\cal T$ is thus a contravariant functor from the category of
closed oriented topological surfaces, morphisms being
covering maps, to the category of finite dimensional complex
manifolds and holomorphic embeddings. We shall have more to say
about this in Section 5.

Over each genus Teichm\"uller space we have a sequence of natural
determinant bundles arising from the powers of the relative
(co-)tangent bundles along the fibers of the universal curve. Indeed,
let ${\om}_g \longrightarrow {\cal C}_g$ be the relative cotangent bundle
for the projection ${\p}_g$ in $(3.2)$. The determinant line bundle
over $\Tg$ arising from its $n$-th tensor power is fundamental,
and we shall denote it by:
$$
DET_{n,g}:= det({\om}^n_{g}) \longrightarrow \Tg, ~~n \in \Z
\eqno(3.5)
$$
Applying Serre duality shows that there is a canonical isomorphism
$DET_{n,g} = DET_{1-n,g}$, for all $n$. $DET_{0,g} = DET_{1,g}$ is
called the {\it Hodge} line bundle over $\Tg$.

These holomorphic line bundles carry natural {\it Quillen hermitian
structure} arising from the Poincar\'e metrics on
the fibers of the universal curve.
Recall that any Riemann surface $Y$ of genus $g\geq 2$ admits a unique
conformal Riemannian metric of constant curvature $-1$, called the
Poincar\'e metric of $Y$. This metric depends smoothly on the conformal
structure, (because of the uniformization theorem with moduli parameters),
and hence, for a family of Riemann surfaces of genus at least
two, the Poincar\'e metric induces a hermitian metric on the relative
tangent/cotangent bundle. We thus obtain Quillen metrics on each
$DET_{n,g}$. The metric functorially assigned by the Quillen metric
on any tensor power of $DET_{n,g}$ will also be referred to as the Quillen
metric on that tensor power.

Observe that by the naturality of the above constructions it follows
that the action of $MC_g$ on $\Tg$ has a natural lifting as unitary
automorphisms of these $DET$ bundles.

We invoke back into play the unramified finite covering
$\pi: \tX \rightarrow X$. Let
$$
{\cal T}(\pi )^*{\cal C}_{\tg}~\longrightarrow ~{\cal T}_g
\eqno{(3.6)}
$$
be the pull-back to ${\cal T}_g$ of the universal family ${\cal C}_{\tg}
\longrightarrow {\cal T}_{\tg}$ using the map ${\cal T}(\pi)$.
Given the topological covering space $\pi$ we therefore
obtain the following \'etale covering map between families of
Riemann surfaces parametrized by $\Tg$:
$$
{\pi}\times id ~:~ {\cal T}(\pi )^*{\cal C}_{\tg}~
\longrightarrow ~{\cal C}_g~:=~X\times {\cal T}_g
$$
This is clearly a holomorphic map. In fact, we have the following
commutative diagram
$$
\matrix{{\cal T}(\pi )^*{\cal C}_{\tg}&
&\mapright{\pi\times id}& & {\cal C}_g\cr &\searrow& &
\swarrow &\cr & &{\cal T}_g& & \cr}
\eqno{(3.7)}
$$
exactly as in the general situation (2.7) above Lemma 2.9.

Now let
$$
id\times {\cal T}(\pi )~:~ {\cal T}(\pi )^*{\cal C}_{\tg}
{}~\longrightarrow ~{\cal C}_{\tg}
$$
denote the tautological lift of the map ${\cal T}(\pi)$.
{}From the definition of the Poincar\'e metric it is clear that
for an unramified covering of Riemann surfaces, $\tY \longrightarrow Y$,
the Poincar\'e metric on $\tY$ is the pull-back of the Poincar\'e
metric on $Y$. If ${\om}_{\tg}$ is the relative cotangent bundle
on ${\cal C}_{\tg}$ then this
compatibility between Poincar\'e metrics implies that the two hermitian
line bundles on ${\cal T}(\pi )^*{\cal C}_{\tg}$ namely,
$(\pi\times id)^*{\om}_g$ and $(id\times{\cal T}(\pi ))^*{\om}_{\tg}$
are canonically isometric.

But since the determinant bundle of a pullback family is the
pullback of the determinant bundle, the holomorphic hermitian bundle
${\cal T}(\pi)^*(det({\om}^n_{\tg})) \longrightarrow {\cal T}_g$ is
canonically isometrically isomorphic to the determinant bundle of
$(id\times {\cal T}(\pi))^*{\om}^n_{\tg} \longrightarrow
{\cal T}(\pi )^*{\cal C}_{\tg}$. Using this and simply applying
Lemma 2.9 to the commutative diagram (3.7) we obtain the following
theorem. (All the Quillen metrics are with respect to the Poincar\'e
metric on fibers.)

\medskip
\noindent{\bf Theorem 3.8a.} {\it The two holomorphic hermitian line
bundles $det({\om}^n_g)^{12.{\rm deg}(\pi )}$ and\\ ${\cal T}(\pi
)^*(det({\om}^n_{\tg}))^{12}$ on ${\cal T}_g$ are canonically isometrically
isomorphic for every integer $n$. In other words, there is a canonical
isometrical line bundle morphism $\Ga(\pi)$ lifting $\Tpi$
and making the following diagram commute:
$$
\matrix{
{\Dng}^{\otimes {12.deg(\pi)}}
&\mapright{{\Ga}(\pi)}
&{\Dntt}
\cr
\mapdown{}
&
&\mapdown{}
\cr
\Tg
&\mapright{\Tpi}
&\Tt
\cr}
$$
}
\medskip

\noindent
{\bf Remark 3.9.} The bundle morphism $\Ga(\pi)$ has been obtained from
Riemann-Roch isomorphisms -- as evinced by the proof of Lemma 2.9. We shall
therefore, in the sequel, refer to these canonical mappings as GRR morphisms.
Tensor powers of the GRR morphisms  will also be referred to as GRR
morphisms. The functoriality of these morphisms is explained below
in Theorem 3.8b.

Let $\bX~\mapright{\r}~\tX~\mapright{\pi}~X$ be two unramified
coverings  between closed surfaces of respective genera
${\bg}$, $\tg$ and $g$. By applying the Teichm\"uller functor we
have the corresponding commuting triangle of embeddings between
the Teichm\"uller spaces:

$$
\matrix{
\Tg
&
&\mapright{\Tpi}
&
&\Tt
\cr
&
\searrow
&
&\swarrow
&
\cr
&&
\Tbg
&&
\cr}
\eqno(3.10)
$$
Here the two slanting embeedings are, of course,
${\cal T}(\pi\circ\r )$ and ${\cal T}(\r)$.
Applying Lemma 2.9(ii) we have

\medskip
\noindent{\bf Theorem 3.8b.}{\it The following triangle of
GRR line bundle morphisms commutes:
$$
\matrix{
\Dng^{\otimes {12.deg(\pi\circ\r)}}
&
&\mapright{}
&
&\Dnt^{\otimes {12.deg(\r)}}
\cr
&
\searrow
&
&\swarrow
&
\cr
&&
\Dnb^{\otimes {12}}
&&
\cr}
\eqno(3.11)
$$
All three maps in the diagram are obtained by applications of Theorem
3.8a, and raising to the appropriate tensor powers.
The triangle above sits over the triangle of Teichm\"uller spaces (3.10),
and the entire triangular prism is a commutative diagram.}

\medskip

\noindent {\bf Remark 3.12.}\,\ The nagging factor of $12$ in Theorems
3.8a and 3.8b can be dealt with as follows. The Teichm\"uller space being a
contractible Stein domain, any two line bundles on ${\cal T}_g$ are
isomorphic. Choose an isomorphism between ${\de} : det({\om}^n_g)^{{\rm
deg}(\pi )} \longrightarrow {\cal T}(\pi )^*(det({\om}^n_{\tg}))$. Hence
$${\de}^{\otimes 12} : det({\om}^n_g)^{{12.\rm deg}(\pi )}
\longrightarrow {\cal T}(\pi )^*(det({\om}^n_{\tg}))^{12}$$ is an
isomorphism. Let $${\tau} : det({\om}^n_g)^{{12.\rm deg}(\pi )}
\longrightarrow {\cal T}(\pi )^*(det({\om}^n_{\tg}))^{12}$$ be the
isomorphism given by the Theorem 3.8a. So $f := {\tau}\circ ({\de}^{\otimes
12})^{-1}$ is a nowhere zero function on ${\cal T}_g$.  Since ${\cal T}_g$
is simply connected, there is a function $h$ on ${\cal T}_g$ such that
$h^{12} = f$. Any two such choices of $h$ will differ by a $12$-th root of
unity. Consider the homomorphism ${\bar \tau} := h.\de$. Clearly ${\bar
\tau}^{\otimes 12} = \tau$. It is easy to see that
for two different choices of the isomorphism $\de$, the two ${\bar\tau}'$s
differ by multiplication with a $12$-th root of unity. Moreover, if we
consider a similar diagram to that in Theorem 3.8b with the factor 12
removed and all the homomorphisms being replaced by the corresponding
analogues of $\bar\tau$, then the diagram commutes up to multiplication
with a $12$-th root of unity.

\noindent{\bf Remark 3.13.}
Recall from above that the action of $MC_g$ in $\Tg$
lifts to the total space of $det({\om}^n_g)$ as bundle
automorphisms preserving the Quillen metric.  There is no action,
a priori, of $MC_g$ on the total space of the the pullback bundle
${\cal T}(\pi )^*(det({\om}^n_{\tg}))$. However, from Theorem 3.8a
the bundle ${\cal T}(\pi )^*(det({\om}^n_{\tg}))^{12}$ gets an
action of $MC_g$ which preserves the pulled back Quillen metric.
Theorem 3.8b ensures the identity between the $MC_g$ actions
obtained by different pullbacks.
\medskip

In [BN] we will consider two special classes of coverings,
namely characteristic covers and cyclic covers. In such situations the map
between Teichm\"uller spaces, induced by the covering, actually descends
to a map between moduli spaces (possibly with level structure).
As mentioned in the Introduction, in that
context we were able to give a proof of the existence of the
GRR morphism of Theorem 3.8a using Weil-Petersson geometry and
topology.

\section{Power law (Principal) bundle morphisms over Teichm\"uller spaces}

We desire to obtain certain canonical geometric objects over the
inductive limits of the finite dimensional Teichm\"uller spaces
by coherently fitting together the determinant line bundles
$DET_{n,g}$ thereon; the limit is taken as $g$ increases by
running through a universal tower of covering maps.
To this end it is necessary to find canonical mappings relating
$\Dng$ to $\Dnt$ where genus $\tg$ covers genus $g$.

Now, given any complex line bundle $\la \ra T$ over any base $T$,
there is a certain canonical mapping of $\la$ to any positive integral
($d$-th) tensor power of itself, given by:
$$
\om_d: \la \longrightarrow {\la}^{\otimes d}
\eqno(4.1)
$$
where $\om_d$ on any fiber of $\la$ is the map $z \mapsto z^{d}$.
Observe that $\om_d$ maps $\la$ minus its zero section to
${\la}^{\otimes d}$ minus its zero section by a map which
is of degree $d$ on the $\CC^{*}$ fibers.
We record the following properties of these maps:

\noindent
{\bf 4.1a.}
The map $\om_d$ is defined independent of any choices of basis,
and it is evidently compatible with base change. [Namely, if we pull
back both $\la$ and $\la^{d}$ over some base $T_1 \rightarrow T$,
then the connecting map $\om_d$ (over $T$) also pulls back to
the corresponding $\om_d$ over $T_1$.]

\noindent
{\bf 4.1b.}
The map $\om_d$ is a homomorphism of the corresponding $\CC^{*}$
principal bundles. When $T$ is a complex manifold, and $\la$ is
a line bundle in that category, then the map $\om_d$ is a holomorphic
morphism between the total spaces of the source and target bundles.

\noindent
{\bf 4.1c.}
If $\la$ is equipped with a hermitian fiber metric, and its tensor
powers are assigned the corresponding hermitian structures, the map
$\om_d$ carries the unit circles to unit circles. (The choice of
a unit circle amongst the natural family of zero-centered
circles in any complex line is clearly equivalent to specifying a
hermitian norm. In this section we will think of hermitian
structure on a line bundle as the choice of a smoothly varying
family of unit circles in the fibers.)

Thus, given a topological covering $\pi:\tX \rightarrow X$,
as in the situation of Theorem 3.8a, we may define a canonical map
$$
\Om(\pi) := \Ga(\pi) \circ \om_{deg(\pi)}:\Dng^{\otimes 12}
\longrightarrow \Dntt
\eqno(4.2)
$$
where $\Ga(\pi)$  is the canonical GRR line bundle morphism
found in Theorem 3.8a. Translating Theorems 3.8a and 3.8b in
terms of these holomorphic maps $\Om(\pi)$ of positive
integral fiber degree, we get:

\medskip
\noindent{\bf Theorem 4.3a.} {\it For each integer $n$, there is
a canonical isometrical holomorphic bundle morphism $\Om(\pi)$
lifting $\Tpi$ and making the following diagram commute:
$$
\matrix{
{\Dngt}
&\mapright{{\Om}(\pi)}
&{\Dntt}
\cr
\mapdown{}
&
&\mapdown{}
\cr
\Tg
&\mapright{\Tpi}
&\Tt
\cr}
$$
By ``isometrical'' we mean that the unit circles of the Quillen
hermitian structures are preserved by the $\Om(\pi)$.}

\smallskip
\noindent{\bf 4.3b.} {\it Let $\pi$ and $\r$ denote two composable
covering spaces between surfaces, as in the situation of Theorem 3.8b.
The following triangle of non-linear isometrical
holomorphic bundle morphisms commutes:
$$
\matrix{
\Dng^{\otimes {12}}
&
&\mapright{}
&
&\Dnt^{\otimes {12}}
\cr
&
\searrow
&
&\swarrow
&
\cr
&&
\Dnb^{\otimes {12}}
&&
\cr}
$$
The horizontal map is $\Om(\pi)$, and the two slanting maps are,
(reading from left to right), $\Om(\pi \circ \r)$ and $\Om(\r)$.
The triangle above sits over the triangle of Teichm\"uller spaces
(3.10), and the entire triangular prism is a commutative diagram.}

The canonical and functorial nature of these connecting maps,
$\Om(\pi)$, will now allow us to produce direct systems of
line/principal bundles over direct systems of Teichm\"uller spaces.

\section{Commensurability Teichm\"uller space and its Automorphism group}

We construct a category $\S$ of certain topological objects and
morphisms: the objects, $Ob(\S)$, are a set of compact oriented topological
surfaces each equipped with a base point ($\star$), there being exactly
one surface of each genus $g \geq 0$; let the object of genus $g$ be
denoted by $X_g$. The morphisms are based isotopy classes of pointed
covering mappings
$$
\pi: (X_\tg, \star) \rightarrow (X_g, \star)
$$
there being one arrow for each such isotopy class. Note that the
monomorphism of fundamental groups induced by (any representative of
the based isotopy class) $\pi$, is unambiguously defined.

Fix a genus $g$ and let $X = X_g$.
Observe that all the morphisms with the fixed target $X_g$:
$$
M_g = \{\a \in Mor(\S): Range(\a)=X_g \}
$$
constitute a {\it directed set} under the partial ordering given by
factorisation of covering maps. Thus if $\a$ and $\b$ are two morphisms from
the above set, then $\b \succ \a$ if and only if the image of the
monomorphism $\pi_1(\b)$ is contained within the image of $\pi_1(\a)$;
that happens if and only if there is a commuting triangle of morphisms of
$\S$ as follows:

$$
\matrix{
X_{g(\b)}
&
&\mapright{\theta}
&
&X_{g(\a)}
\cr
&
\searrow
&
&\swarrow
&
\cr
&&
X_g
&&
\cr}
$$
Here $X_{g(\a)}$ denotes the domain surface for $\a$ (similarly
$X_{g(\b)}$), and the two slanting arrows are (reading from left to
right), $\b$ and $\a$. It is important to note that the factoring
morphism $\theta$ is {\it uniquely} determined because we are working
with base points. The directed property of $M_g$ follows by a simple
fiber-product argument.
[Remark: Notice that the object of genus $1$ in $\S$ only has
morphisms to itself -- so that this object together with all its
morphisms (to and from) form a subcategory.]

Recall from Section 3 that each morphism of $\S$ induces a
proper, holomorphic, Teichm\"uller-metric preserving
embedding between the corresponding finite-dimensional Teichm\"uller
spaces. We can thus create the natural {\it direct system of Teichm\"uller
spaces} over the above directed set $M_g$, by associating to each $\a
\in M_g$ the Teichm\"uller space $\T(X_{g(\a)})$, and for each $\b
\succ \a$ the corresponding holomorphic embedding $\T(\theta)$ (with
$\theta$ as in the diagram above).
Consequently, we may form the {\it direct limit Teichm\"uller space over
$X=X_g$}:
$$
\Tin(X_g) = \TinX := ind. lim. \T(X_{g(\a)})
\eqno(5.1)
$$
the inductive limit being taken over all $\a$ in the directed set
$M_g$. This is our {\it commensurability Teichm\"uller space}.

\noindent
{\bf Remark:} Over the same directed set $M_g$ we may also define a
natural {\it inverse system of surfaces}, by asscoiating to $\a \in
M_g$ a certain copy, $S_{\a}$ of the pointed surface $X_{g(\a)}$.
[Fix a universal covering over $X=X_g$. $S_{\a}$ can be taken to be the
universal covering quotiented by the action of the subgroup
$Im(\pi_1(\a)) \subset {\pi_1}(X,\star)$.] If $g \ge 2$, then
the inverse limit of this system is the {\it universal solenoidal surface}
$H_{\infty}$ whose Teichm\"uller theory was studied in [S], [NS].
The completion of $\TinX$ in the Teichm\"uller metric is
$\T(H_{\infty})$.

A remarkable but obvious fact about this construction is that every
morphism $\pi:Y \rightarrow X$ of $\S$ induces a natural Teichm\"uller
metric preserving {\it homeomorphism}
$$
\Tin(\pi): \TinY \longrightarrow \TinX
\eqno(5.2)
$$
$\Tin(\pi)$ is invertible simply because the morphisms
of $\S$ with target $Y$ are cofinal with those having target $X$.
If we consider objects and maps to be
continuous/holomorphic on the inductive limit spaces when they are
continuous/holomorphic when restricted to the finite-dimensional
strata, then it is clear that $\Tin(\pi)$ is a biholomorphic
identification. (Note that $\Tin$ acts covariantly, since it is defined by
a morphism of direct systems, although the Teichm\"uller functor $\T$
of (3.3) was contravariant.)

It follows that each $\TinX$, (and its completion $\T(H_{\infty})$),
is equipped with a large {\it automorphism group} -- one from each
(undirected) cycle of morphisms of $\S$ starting from $X$ and
returning to $X$. By repeatedly using pull-back diagrams (i.e., by
choosing  the appropriate connected component of the fiber product of
covering maps), it is easy to see that the automorphism arising from
any (many arrows) cycle can be obtained simply from a two-arrow cycle
$\tX {\ra \atop \ra} X$. Namely, whenever we have
(the isotopy class of) a ``self-correspondence'' of $X$ given by
two non-isotopic coverings, say $\a$ and $\b$,
$$
\tX {\ra \atop \ra} X
\eqno(5.3)
$$
we can create an automorphism of $\TinX$ defined as the composition:
${\Tin(\b)}\circ{(\Tin(\a))^{-1}}$.  Therefore each of these automorphisms
-- arising from any arbitrarily complicated cycle of coverings (starting
and ending at $X$) -- is obtained as one of these simple
``two-arrow'' compositions. These automorphisms form a group that we
shall call the {\it commensurability modular group}, $CM_{\infty}(X)$,
acting on the universal commensurability Teichm\"uller space $\TinX$.

We make some further remarks regarding this large new mapping class
group. Consider the abstract graph ($1$-complex), $\Ga(\S)$,
obtained from the category $\S$ by looking at the objects
as vertices and the (undirected) arrows as edges. It is
clear from the definition above that the fundamental group
of this graph, viz. $\pi_{1}(\Ga(\S),X)$, is acting on $\TinX$
as these automorphisms. In fact, we may fill in all the ``commuting
triangles''  -- i.e., fill in the $2$-cells in this abstract
graph whenever two morphisms (edges) compose to give a third edge;
the thereby-reduced fundamental group of this $2$-complex produces
on $\TinX$ the action of $CM_{\infty}(X)$.

It is interesting to observe that this new modular group $CM_{\ify}(X)$
of automorphisms on $\TinX$ corresponds exactly to ``virtual automorphisms''
of the fundamental group $\pi_{1}(X)$, -- generalizing the
classical situation where the usual automorphism group $Aut(\pi_{1}(X))$
appears as the action via modular automorphisms on $\T(X)$.

Indeed, given any group $G$, one may define its associated group of
``virtual'' automorphisms; as opposed to usual automorphisms,
for virtual automorphisms we demand that they be defined only
on some finite index subgroup of $G$. To be precise,
consider isomorphisms $\r:H \ra K$ where $H$ and $K$ are
subgroups of finite index in $G$. Two such isomorphisms (say $\r_1$
and $\r_2$) are considered equivalent if there is a finite index subgroup
(sitting in the intersection of the two domain groups) on which they
coincide. The equivalence class $[\r]$  -- which is like the {\it
germ} of the isomorphism $\r$ -- is called a {\it virtual automorphism}
of $G$; clearly the virtual automorphisms of $G$ constitute a group,
$Vaut(G)$, under the obvious law of composition, (namely, compose
after passing to deeper finite index subgroups, if necessary).

We shall apply this concept to the fundamental group of a surface of
genus $g$, ($g>1$). It is clear from definition that the group
$Vaut(\pi_{1}(X_g))$ {\it is genus independent}, as is to be expected
in our constructions.

In fact, $Vaut$ presents us a neat way of formalizing the ``two-arrow
cycles'' which we introduced to represent elements of $\CM$.
Letting $G = \pi_{1}(X)$, (recall that $X$ is already equipped with a
base point), the two-arrow diagram (5.3) above corresponds to the
following well-defined virtual automorphism of $G$:
$$
[\r] = [{\b}_{*}\circ{\a}_{*}^{-1}:{\a}_{*}(\pi_{1}(\tX)) \rightarrow
{\b}_{*}(\pi_{1}(\tX))]
$$
Here ${\a}_{*}$ denotes the monomorphism of the fundamental group
$\pi_{1}(\tX)$ into $\pi_{1}(X) = G$, and similarly ${\b}_{*}$.
We let $Vaut^{+}({\pi}_{1}(X))$ denote the subgroup of $Vaut$ arising
from pairs of orientation preserving coverings. The final upshot is
that  $\CMX$ is {\it isomorphic} to $Vaut^{+}(\pi_{1}(X))$ and there
is a natural surjective homomorphism:
${\pi}_{1}(\Ga(\S),X) \rightarrow Vaut^{+}({\pi}_{1}(X))$ whose kernel is
obtained by filling in all commuting triangles in $\Ga(\S)$.

\noindent
{\it Acknowledgement:} The concept of $Vaut$ has arisen in group
theory papers -- for example [Ma],[MT]. We are grateful to Chris Odden for
pointing out these references.

\noindent {\bf Remark 5.4.} For the genus one object $X_1$ in $\S$,
we know that the Teichm\"uller spaces for all unramified coverings
are each a copy of the upper half-plane $H$. The maps $\Tpi$ are
M\"obius identifications of copies of the half-plane with itself,
and we easily see that the pair $(\Tin(X_1),CM_{\ify}(X_1))$ is
identifiable as $(H,PGL(2,\QQ))$. In fact, $GL(2,\QQ) \cong Vaut(\Z \oplus
\Z)$,
and $Vaut^{+}$ is precisely the subgroup of index $2$ therein, as expected.
Notice that the action has dense orbits in the genus one case.

On the other hand, if $X \in Ob(\S)$ is of any genus $g \geq 2$,
then we get an infinite dimensional ``ind-space'' as $\TinX$ with
the action of $\G(X)$ on it as described. Since the tower of coverings
over $X$ and $Y$ (both of genus higher than $1$) eventually become
cofinal, it is clear that {\it for any choice of genus higher than one we get
one isomorphism class of pairs} $(\Tin, CM_{\ify})$. (It is not known
whether the action has dense orbits in this situation; this matter is
related to some old queries on coverings of Riemann surfaces.)

We work now over the direct system of the higher genus example
$(T_{\ify}, CM_{\ify})$ and obtain the main theorem. We will first
explain some preliminary material on direct limits of holomorphic
line bundles over a direct system of complex manifolds.

Given a direct system $T_{\a}$ of complex manifolds, and line
bundles $\xi_{\a}$ over these, suppose that there are power
law maps as the $\Om(\pi)$ above, between the corresponding principal
$\C*$ bundles covering the mappings in the direct system of
base manifolds.

Let $N$ denote the directed set of positive integers ordered by
divisibility. For each $\la \in N$ take a copy of $\C*$, call it
$(\C*,\la)$ and form the direct system $\{ (\C* ,\la)\}$ where
$(\C*,\la) \ra (\C* ,{\la}')$ is given by the power law map:
$z \ra z^d$ when ${\la}'=d{\la}$. These maps are homomorphisms
of groups, and the direct limit over $N$ is canonically isomorphic to
the group $\CQ := \CC\otimes_{\Z}\QQ$. [The isomorphism maps the
equivalence class of the element $(z,\la) \in (\C*,\la)$
to $z \otimes {1/\la} \in \CQ$.]
The direct limit object obtained from the power law connecting maps
between the principal bundles associated to the ${DET}_{n}^{12}$ system
over the Teichm\"uller spaces will give us a $\CQ$ principal bundle
over the universal commensurability Teichm\"uller space $\Tin$, at least at
the level of sets. The topological and holomorphic structure on these
sets is defined for maps into these objects which factor through
the direct system by imposing these properties on the factorizations.

Let us consider the direct limit bundles obtained from a family of such
bundles $\xi_{\a}$, and from the family obtained by raising each
$\xi_{\a}$ to the tensor power $d$. These are two $\CQ$ bundles over
the direct limit of the bases which may be thought to have the
same total spaces (as sets) but the $\CQ$ action on the second one is
obtained by precomposing the original action by the automorphism
of $\CQ$ obtained from the homomorphism $z \mapsto z^d$ on $\C*$.

\noindent
{\bf Theorem 5.5.} {\it Fix any integer $n$. Starting from any
base surface $X \in Ob(\S)$, we obtain a direct system of principal
$\CC^*$ bundles $\L_n(Y) := DET_{n,g(Y)}^{\otimes 12}$ over the
Teichm\"uller spaces $\T(Y)$ with holomorphic homomorphisms $\Om(\pi)$
(see Theorem 4.3) between the total spaces; here $Y \mapright{\pi} X$
is an arbitrary morphism of $\S$ with target $X$.

By passing to the direct limit, one
therefore obtains over the universal commensurability Teichm\"uller space,
$\TinX$, a principal $\CQ$ bundle:
$$
\L_{n,\infty}(X) = ind. lim. \L_n(Y)
$$
Since the maps $\Om(\pi)$ preserved the Quillen unit circles, the limit
object also inherits such a Quillen ``hermitian'' structure.

The construction is functorial with respect to change of the base $X$
in the obvious sense that the directed systems and their limits are
compatible with the biholomorphic identifications $\Tin(\pi)$ of
equation (5.2). In particular, the commensurability modular group
action $CM_{\ify}(X)$ on $\TinX$ has a natural lifting to
$\L_{n,\infty}(X)$ -- acting by unitary automorphisms.

Finally, the Mumford isomorphisms persist:
$$
\L_{n,\infty}(X) = \L_{0,\infty}(X)^{\otimes (6n^{2} - 6n + 1)}
$$
Namely, if we change the action of $\CQ$ on the ``Hodge'' bundle
$\L_{0,\infty}$ by the ``raising to the $(6n^2-6n+1)$ power''
automorphism of $\CQ$, then the principal $\CQ$ bundles are
canonically isomorphic.
}

\noindent {\bf Remark 5.6.}
In other words,
the Mumford isomorphism in the above theorem means that
$\L_{n,\infty}$ and $\L_{0,\infty}$ are equivariantly isomorphic
relative to the automorphism of $\CQ$ induced by the homomorphism
of $\C*$ that raises to the power exhibited.
Also, we could have used the Quillen hermitian structure to reduce the
structure group from $\CC^{*}$ to $U(1)$, and thus obtain direct systems
of $U(1)$ bundles over the Teichm\"uller spaces. Passing to the direct
limit would then produce $U(1) \otimes \QQ :=$ {\it ``tiny circle''}
bundles over $\Tin$, which can be tested for maps into these objects
as above.

\smallskip
\noindent
{\bf Rational line bundles on ind-spaces:} A line bundle on the inductive
limit of an inductive system of varieties or spaces, is, by definition
([Sha]), a collection of line bundles on each stratum (i.e., each member
of the inductive system of spaces) together with compatible line bundle
(linear on fibers) morphisms. Such a direct system of line bundles
determines an element of the inverse limit of the Picard groups of
the stratifying spaces. See [KNR], [Sha].
(Recall: For any complex space $T$, $Pic(T)$ := the group (under $\otimes$)
of isomorphism classes of line bundles on $T$. In the case of the
Teichm\"uller spaces, we refer to the modular-group invariant
bundles as constituting the relevant Picard group -- see [BN].)

Now, utilising the GRR morphisms $\Ga(\pi)$ themselves, (without
involving the power law maps), we know from Section 3 that the
``$d$-th root'' of the bundle $\Dnt$ fits together with the bundle
$\Dng$ ($d=(\tg-1)/(g-1)$).
A ``rational'' line bundle over the inductive limit is defined to be
an element of the inverse limit of the $Pic_{\QQ} = Pic \otimes {\QQ}$.
Therefore we may also state a result about the existence of
canonical elements of the inverse limit,
$\lim_{\leftarrow}Pic({\cal T}_{g_i})_{\QQ}$, by our construction.
Indeed, in the notation of Section 3, by using the morphisms
${\Ga(\pi)}\ot{1/{deg(\pi)}}$ between $\Dng$ and ${\Dnt}\ot{1/{deg(\pi)}}$
to create a directed system, we obtain canonical elements representing
the Hodge and higher $DET_n$ bundles, with respective Quillen metrics:
$$ {\La}_{m} \in
\lim_{\leftarrow}Pic({\cal T}_{g_i})_{\QQ}, ~~m=0,1,2,..
\eqno(5.7)
$$
The pullback (i.e., restriction) of $\La_m$ to each of the stratifying
Teichm\"uller spaces ${\cal T}_{g_i}$ is $(n_i)^{-1}$ times the
corresponding $DET_{m}$ bundle, (with $(n_i)^{-1}$ times its Quillen metric),
over ${\cal T}_{g_i}$. Here $n_i$ is the degree of the covering of the
surface of genus $g_i$ over the base surface. As rational hermitian line
bundles the Mumford isomorphisms persist:
$$
\La_{m} = {\La_{0}}^{\otimes (6m^{2} - 6m +1)}
\eqno(5.8)
$$
as desired. This statement is different from that of the Theorem.
For further details see [BN].

\smallskip

\noindent{\bf Polyakov measure on $\Mg$ and our constructions:} In his
study of bosonic string theory, Polyakov constructed a measure
on the moduli space ${\cal M}_g$ of curves of genus $g(\geq 2)$.
Details can found, for example, in
[Alv], [Nel]. Subsequently, Belavin and Knizhnik, [BK],
showed that the Polyakov measure has the following elegant mathematical
description. First note that a hermitian metric on the canonical bundle of
a complex space gives a measure on that space. Fixing a volume form
(up to scale) on a space therefore amounts to fixing a fiber
metric (up to scale) on the canonical line bundle, $K$, over
that space. But the Hodge bundle $\l$ has its natural Hodge metric
(arising from the $L^2$ pairing of holomorphic
1-forms on Riemann surfaces). Therefore we may transport the corresponding
metric on ${\l}^{13}$ to $K$ by Mumford's isomorphism, (as we know the
choice of this isomorphism is unique up to scalar) -- thereby obtaining a
volume form on ${\cal{M}}_g$. [BK] showed that this is none other than the
Polyakov volume.  Therefore, the presence of Mumford isomorphisms over the
moduli space of genus $g$ Riemann surfaces describes the Polyakov measure
structure thereon.

Above we have succeeded in fitting together the Hodge and higher $DET$
bundles over the ind space $\Tin$, together with the relating Mumford
isomorphisms. We thus have from our results a structure on $\Tin$ that
suggests a genus-independent, universal, version of the Polyakov
structure.

We remark that since the genus is considered the
perturbation parameter in the above formulation of the standard
perturbative bosonic Polyakov string theory, our work can be considered as
a contribution towards a {\it non-perturbative} formulation of that theory.

\newpage


\noindent
Authors' Addresses: \\
\noindent
Tata Institute of Fundamental Research,\\
Colaba, Bombay 400 005, INDIA; ``indranil@math.tifr.res.in''

\noindent
The Institute of Mathematical Sciences,\\
CIT Campus, Madras 600 113, INDIA; ``nag@imsc.ernet.in''

\noindent
Institut des Hautes Etudes Scientifiques,\\
91440 Bures-sur-Yvette, FRANCE; ``sullivan@ihes.fr''


\begin{thebibliography}{99}

\bibitem[Alv]{Alv} O. Alvarez, Theory of strings with boundaries
: fluctuations, topology and quantum geometry, {\it Nuc. Phys.
B216}, (1983), 125-184.

\bibitem[BGS]{BGS} J.Bismut, H.Gillet and C.Soul\'e, Analytic
torsion and holomorphic determinant bundles, I,II,III, {\it
Comm. Math Phys., 115}, 1988, 49-126, 301-351.

\bibitem[BK]{BK} A. Belavin, V. Knizhnik, Complex geometry and
quantum string theory. {\it Phys. Lett., 168 B}, (1986),
201-206.

\bibitem[BN]{BN} I. Biswas and S. Nag, Weil-Petersson geometry
and determinant bundles over inductive limits of moduli spaces,
{\it preprint}.

\bibitem[Bos]{Bos} J.B. Bost, Fibres determinants, determinants
regularises et mesures sur les espaces de modules des courbes
complexes, {\it Semin.  Bourbaki, 152-153}, (1987), 113-149.

\bibitem[D]{D} P. Deligne, Le d\'eterminant de la cohomologie,
{\it Contemporary Math., 67}, (1987), 93-177.

\bibitem[GH]{GH} P. Griffiths, J. Harris, Principles of
algebraic geometry, Wiley-Interscience, John Wiley
\& Sons, New York.

\bibitem[GS]{GS} H. Gillet and C. Soul\'e, An arithmetic
Riemann-Roch theorem, {\it Inv. Math., 110}, (1992), 473-543.

\bibitem[KM]{KM} F. Knudsen and D. Mumford, The projectivity of
the moduli space of stable curves I, preliminaries on ``det''
and ``div'', {\it Math. Scand. 39},
(1976) 19-55.

\bibitem[KNR]{KNR} S. Kumar, M.S. Narasimhan, A. Ramanathan,
Infinite Grassmannians and moduli spaces of $G$-bundles, {\it
Math. Ann. 300}, (1994), 41-73.

\bibitem[Ma]{Ma} A. Mann, Problem about automorphisms of infinite groups,
{\it Second Intn'l Conf Group Theory}, Debrecen, (1987).

\bibitem[MT]{MT} F. Menegazzo and M. Tomkinson, Groups with trivial
virtual automorphism group, {\it Israel Jour. Math, 71}, (1990), 297-308.

\bibitem[M]{M} L. Moret-Bailly, La formule de Noether pour les
surfaces arithm\'etiques, {\it Invent. Math., 98}, (1989),
491-498.

\bibitem[Mum]{Mum} D. Mumford, Stability of projective varieties,
{\it L'Enseign. Math., 23}, (1977), 39-100.

\bibitem[N]{N} S. Nag, The Complex Analytic Theory of Teichm\"uller
Spaces, Wiley-Interscience, John Wiley \& Sons, New York, 1988.

\bibitem[NS]{NS} S. Nag, D. Sullivan, Teichm\"uller theory and
the universal period mapping via quantum calculus and the
$H^{1/2}$ space on the circle, {\it Osaka J. Math., 32, no.1},
(1995), 1-34.  [Max-Planck-Inst. Bonn preprint MPI-94-54.]

\bibitem[Nel]{Nel} P.Nelson, Lectures on strings and moduli
space, {\it Physics Reports}, 1987.

\bibitem[Q]{Q} D. Quillen, Determinants of Cauchy-Riemann
operators over a Riemann surface. {\it Func. Anal. Appl., 19},
(1985), 31-34.

\bibitem[S]{S} D. Sullivan, Relating the universalities of
Milnor-Thurston, Feigenbaum and Ahlfors-Bers, in Milnor
Festschrift {\it Topological Methods in Modern mathematics} (ed.
L. Goldberg, A. Phillips) Publish or Perish, (1993), 543-563.

\bibitem[SGA]{SGA} Seminaire de Geometrie Algebrique, 1.,
Rev\'etments \'etales et groupe fondamental, Lecture Notes in Math.,
no. 224, Springer-Verlag, Berlin.

\bibitem[Sha]{Sha} I.R. Shafarevich, On some
infinite-dimensional groups II, {\it Math USSR Izvest., 18},
(1982), 185-194.

\bibitem[Wol]{Wol} S. Wolpert, Chern forms and the Riemann
tensor for the moduli space of curves, {\it Invent. Math., 85},
(1986), 119-145.

\bibitem[ZT]{ZT} P.G. Zograf and L.A. Takhtadzhyan, A local
index theorem for families of {$\bar \partial$}- operators on
Riemann surfaces, {\it Russian Math Surveys, 42}, (1987),
169-190.

\end{thebibliography}
\end{document}
